\input amstex
\documentstyle{amams} 
\annalsline{156}{2002}
\received{July 3, 2001}
\startingpage{333}
\catcode`\@=11
\font\twelvemsb=msbm10 scaled 1100

\font\ninemsb=msbm10 scaled 800
\newfam\msbfam
\textfont\msbfam=\twelvemsb  \scriptfont\msbfam=\ninemsb
  \scriptscriptfont\msbfam=\ninemsb
\def\msb@{\hexnumber@\msbfam}
\def\Bbb{\relax\ifmmode\let\next\Bbb@\else
 \def\next{\errmessage{Use \string\Bbb\space only in math
mode}}\fi\next}
\def\Bbb@#1{{\Bbb@@{#1}}}
\def\Bbb@@#1{\fam\msbfam#1}
\catcode`\@=12

 \catcode`\@=11
\font\twelveeuf=eufm10 scaled 1100
\font\teneuf=eufm10
\font\nineeuf=eufm7 scaled 1100
\newfam\euffam
\textfont\euffam=\twelveeuf  \scriptfont\euffam=\teneuf
  \scriptscriptfont\euffam=\nineeuf
\def\euf@{\hexnumber@\euffam}
\def\frak{\relax\ifmmode\let\next\frak@\else
 \def\next{\errmessage{Use \string\frak\space only in math
mode}}\fi\next}
\def\frak@#1{{\frak@@{#1}}}
\def\frak@@#1{\fam\euffam#1}
\catcode`\@=12

\def\Q{\ifx\MYUN\Bbb {\font\bb=msbm10 \hbox{\bb Q}} \else {\Bbb Q} \fi}
 \def\irr#1{{\rm Irr}(#1)}
\def\cent#1#2{{\bf C}_{#1}(#2)}
\def\syl#1#2{{\rm Syl}_#1(#2)}
\def\gen#1{\langle#1\rangle}
\def\norm#1#2{{\bf N}_{#1}(#2)}

\def\sbs{\subseteq}
\def\ref#1{{[#1]}}

\def\alperin{1}
\def\broue{2}
\def\dade{3}
\def\fongsrin{4}
\def\fong{5}
\def\brown{6}
\def\book{7}
\def\isnav{8}
\def\mckay{9}
\def\mckaya{10}
\def\wilson{11}

\def\cl#1#2{{\rm cl}_{#1}(#2)}
\def\irrpp#1{{\rm Irr}_{p'}(#1)}

\title{New refinements of the McKay conjecture\\ for arbitrary finite groups} 
\shorttitle{New refinements of the McKay conjecture}  
  \twoauthors{I.~M.~Isaacs}{Gabriel Navarro}
\institutions{University of Wisconsin, Madison, WI\\
{\eightpoint {\it E-mail address\/}: isaacs@math.wisc.edu}\\
\vglue6pt
Universitat de Valencia,  Valencia, Spain\\
{\eightpoint {\it E-mail address\/}:  gabriel@uv.es}}

\centerline{\bf Abstract}
\vglue12pt

Let $G$ be an arbitrary finite group and fix a prime number $p$. The
McKay conjecture asserts that $G$ and the normalizer in $G$ of a Sylow
$p$-subgroup have equal numbers of irreducible characters with degrees
not divisible by $p$. The Alperin-McKay conjecture is version of this
as applied to individual Brauer $p$-blocks of $G$. We offer evidence
that perhaps much stronger forms of both of these conjectures are true.
 
\vglue-3pt
\section{Introduction and Conjecture A}

Let $G$ be an arbitrary finite group and fix a prime number $p$. As is well
known, there seem to be some mysterious and unexplained connections
between the representation theory of $G$ and that of certain of its
$p$-local subgroups. For example, it appears to be true that if $P$ is a
Sylow $p$-subgroup of $G$ and $N = \norm GP$, then equal numbers of the
irreducible (complex) characters of $G$ and of $N$ have degrees not
divisible by $p$. This ``McKay Conjecture" was first proposed by J.~McKay
in \ref\mckay, where it was stated in the case where $G$ is simple and
$p = 2$. (See also \ref\mckaya.) The more general formulation of the
conjecture, for arbitrary finite groups and arbitrary primes, was stated by
J.~L.~Alperin in \ref\alperin, although it was first suggested by the first
author in \ref\brown, where it was proved for (solvable) groups of odd
order. (In fact, in the odd-order case considered in \ref\brown, a natural
bijection was constructed between the sets of irreducible characters of
$p'$-degree of $G$ and of $\norm GP$, but it is known that no natural
bijection exists in general.) The McKay conjecture has been verified for
several additional types of groups including $p$-solvable groups, symmetric
groups and the sporadic simple groups. As yet, however, no one has given a
proof, or has even proposed an explanation for why this conjecture might
hold in the general case.

One could argue that the more precisely a conjecture can be stated, the
better it will be understood and thus, perhaps, the easier it might become
to discover a proof. In fact, the McKay conjecture was generalized and
strengthened by Alperin, who formulated   a version that applies to the
Brauer $p$-blocks of $G$. (The Alperin-McKay conjecture was first
proposed in \ref\alperin. We will review its statement in Section 2, where
we discuss blocks.)

In this note we propose several further refinements of the McKay
conjecture and the Alperin-McKay conjecture. To state the first of these,
we define the integer $M_k(G)$, where $G$ is an arbitrary finite group and
$k$ is an integer not divisible by our fixed prime $p$. We write $M_k(G)$ to
denote the total number of irreducible characters of $G$ having degree
congruent modulo $p$ to $\pm k$. For example, if $p = 5$ and $G$ is the
alternating group $A_5$, we see that $M_1(G) = 2$ since $A_5$ has one
irreducible character of degree $1$ and one of degree $4$. Also,
$M_2(G) = 2$ since $A_5$ has two irreducible characters of degree~$3$.
(Note that for odd primes~$p$, the only values of $k$ that we need to
consider are $1 \le k \le (p-1)/2$.)

\nonumproclaim{Conjecture A}
 Let $G$ be an arbitrary finite group and let
$N = \norm GP${\rm ,} where $P \in \syl pG$. Then for each integer $k$ not
divisible by $p${\rm ,} we have $M_k(G) = M_k(N)$.
\endproclaim

For example, if $p = 5$ and $G = A_5$, we saw that $M_1(G) = 2 = M_2(G)$.
In this case, we know that $N$ is dihedral of order $10$, and thus $N$ has
two irreducible characters of degree $1$ and two of degree $2$ and we see
that $M_1(N) = 2 = M_2(N)$, as predicted by Conjecture~A.

Note that if $p = 2$ or $p = 3$, then $M_1(G)$ is the number of irreducible
characters of $G$ of degree not divisible by $p$, and so for these two
primes, Conjecture~A is exactly equivalent to the McKay conjecture. For
$p > 3$, however, we see that the number of irreducible characters of $G$
of degree not divisible by $p$ is the sum of the numbers $M_k(G)$ for
$1 \le k \le (p-1)/2$, and so for these primes, Conjecture~A is strictly
stronger than the McKay conjecture.

But what is the evidence that our conjecture is true?  In the case where 
$G$ has odd order, the bijection constructed in \ref\brown\ carries a
$p'$-degree character $\chi \in \irr G$ to a character $\xi \in \irr{\norm GP}$
such that $\chi(1) \equiv \pm \xi(1)$ mod $p$, and hence Conjecture~A
certainly holds in this case. More generally, but using some deep theory,
the conjecture can be proved for all $p$-solvable groups. Also, P.~Fong
\ref\fong\ has recently succeeded in proving it for all symmetric groups (for
all primes). Furthermore, Conjecture~A holds if the Sylow $p$-subgroup is
cyclic. (This follows from E.~C.~Dade's cyclic defect theory \ref\dade, and
we shall have more to say about that in Section~2.)

If $G = {\rm GL}(n,q)$, where $q$ is a power of the prime $p$, then it is known
that all irreducible characters of $G$ of degree not divisible by $p$ have
degrees congruent to $\pm1$ mod $p$. (In fact, these degrees are
congruent to $\pm1$ mod $q$; this fact follows from results in
\ref\fongsrin.) Also, it is not too hard to see that all of the irreducible
$p'$-degrees for the relevant Sylow normalizer of this group are powers of
$q - 1$, and hence they too are all congruent to $\pm1$ mod $p$. Since the
McKay conjecture is known to be true for ${\rm GL}(n,q)$ in the defining
characteristic, (see \ref\alperin), it follows that Conjecture~A is also true
in this case. (In fact, with some additional work, the conjecture can also be
checked for ${\rm SL}(n,q)$ in the defining characteristic.)

Finally, we mention that Conjecture~A is also true for all primes for all of
the sporadic simple groups. Since the McKay conjecture is known to hold
for these groups \ref\wilson\ and since Conjecture~A is automatically true
when the Sylow $p$-subgroup is cyclic, we see that it suffices to check the
conjecture for primes exceeding $3$ and for which a Sylow subgroup is not
cyclic. We have carried out this check, relying on the ATLAS for the
irreducible character degrees and the paper \ref\wilson\ of R.~A.~Wilson
for the Sylow normalizers of these groups. (However, Wilson's paper has an
error: the normalizer of a Sylow $5$-subgroup in $Fi_{23}$ is incorrectly
described. When we reported this to Wilson, he provided us with a
corrected version.)

The following table gives the relevant data. The third column lists the
numbers $M_k(G)$ for $1 \le k \le (p-1)/2$. The McKay conjecture, of
course, predicts only the sum of these numbers, while Conjecture~A
predicts each of the numbers in the third column.
\medskip
\centerline{
\vtop{\offinterlineskip
\hrule
\halign{&\strut\vrule~\hfil#\hfil~&\vrule~\hfil#\hfil~&\vrule%
~\hfil#~&~\hfil#~&~~\hfil#~\hfil\vrule\cr
Group&Prime&&&\cr
\noalign{\hrule}
$J_2$&5&12&2&\cr
\noalign{\hrule}
$HS$&5&~9&4&\cr
\noalign{\hrule}
$McL$&5&~9&4&\cr
\noalign{\hrule}
$He$&5&~8&8&\cr
\noalign{\hrule}
&7&12&7&1\cr
\noalign{\hrule}
$Ru$&5&10&10&\cr
\noalign{\hrule}
$Suz$&5&~8&8&\cr
\noalign{\hrule}
$O$'$N$&7&12&7&1\cr
\noalign{\hrule}
$CO_3$&5&10&10&\cr
\noalign{\hrule}
$CO_2$&5&10&10&\cr
\noalign{\hrule}
$Fi_{22}$&5&10&10&\cr
\noalign{\hrule}
$HN$&5&10&10&\cr
\noalign{\hrule}
$Ly$&5&25&0&\cr
\noalign{\hrule}
}}\quad\quad
\vtop{\offinterlineskip
\hrule
\halign{&\strut\vrule~\hfil#\hfil~&\vrule~\hfil#\hfil~&\vrule%
~\hfil#~&~\hfil#~&~\hfil#~&~\hfil#~&~\hfil#~&~~\hfil#~\hfil\vrule\cr
Group&Prime&&&&&&\cr
\noalign{\hrule}
$Th$&5&10&10&&&&\cr
\noalign{\hrule}
&7&9&9&9&&&\cr
\noalign{\hrule}
$Fi_{23}$&5&20&20&&&&\cr
\noalign{\hrule}
$Co_1$&5&25&0&&&&\cr
\noalign{\hrule}
&7&9&9&9&&&\cr
\noalign{\hrule}
$J_4$&11&12&15&10&5&0&\cr
\noalign{\hrule}
$Fi_{24}'$&5&28&28&&&&\cr
\noalign{\hrule}
&7&12&9&4&&&\cr
\noalign{\hrule}
$B$&5&25&0&&&&\cr
\noalign{\hrule}
&7&27&27&27&&&\cr
\noalign{\hrule}
$M$&5&40&40&&&&\cr
\noalign{\hrule}
&7&49&0&0&&&\cr
\noalign{\hrule}
&11&10&10&10&10&10&\cr
\noalign{\hrule}
&13&12&18&12&6&3&4\cr
\noalign{\hrule}
}}}
\vglue18pt

Conjecture~A has an amusing consequence for symmetric groups. (As we
mentioned, this case of our conjecture has been proved by Fong.) The fact
is that if $n \ge p$, then all of the numbers $M_k(S_n)$ are multiples of the
prime $p$. To prove this, of course, it suffices to check that the numbers
$M_k(N)$ are multiples of $p$, where $N$ is the normalizer of a Sylow
$p$-subgroup $P$ of the symmetric group $S_n$. It is not very hard to
establish this fact using known information about the Sylow normalizers in
symmetric groups. (See Fong's paper \ref\fong\ for more detail.)

Finally, we mention that there is a significant difference between
Conjecture~A and the many other conjectures that relate the
representation theory of an arbitrary finite group to the $p$-local structure
of the group. (In addition to the McKay conjecture, these include the
Brauer height conjecture, the Alperin weight conjecture and the several
variations and strengthenings of the latter that were formulated by
E.~C.~Dade.) All of these conjectures deal only with the $p$-parts of
irreducible character degrees while our Conjecture~A, of course, is
concerned also with $p'$-parts of character degrees. In the following
section, we discuss Conjecture~B, which to an even greater extent is also
concerned with $p'$-parts of character degrees. 

\section{Blocks}

Let $B$ be a $p$-block of an arbitrary finite group $G$ and let $D$ be a
defect group for $B$. (Recall that $D$ is a $p$-subgroup of $G$ that is
uniquely determined up to $G$-conjugacy by $B$.) As is customary, we will
write $\irr B$ to denote the subset of $\irr G$ consisting of those characters
that belong to the block $B$. Recall that the degrees of the characters in
$\irr B$ are all divisible by $|P|/|D|$, where $P$ is a Sylow $p$-subgroup
of $G$, and that those members of $\irr B$ whose $p$-part is exactly equal
to $|P|/|D|$ are the ``height zero" characters of $B$. 

Now let $N = \norm GD$. Recall that R.~Brauer's famous First Main Theorem
asserts that there is a certain natural bijection between the set of
$p$-blocks of $G$ having defect group $D$ and the set of $p$-blocks of $N$
having defect group $D$. If the block $b$ of $N$ corresponds to the block
$B$ of $G$ under this bijection, we say that $b$ is  the ``Brauer
correspondent" of $B$ with respect to the defect group~$D$. 

Consider the case where $D \in \syl pG$, so that $B$ and $b$ are blocks of
``maximal defect". It is easy to see that the members of $\irr G$ having
degree not divisible by $p$ are exactly all of the height zero characters in
all $p$-blocks of $G$ that have maximal defect. Similarly, the members of
$\irr N$ with degree not divisible by $p$ are just the height zero characters
in the maximal defect $p$-blocks of $N$. The McKay conjecture asserts,
therefore, that the total number of height zero characters in $p$-blocks of
maximal defect of $G$ is equal to the total number of height zero
characters in $p$-blocks of maximal defect of $N$. Since the latter blocks
are exactly the Brauer correspondents of the former, it is reasonable to
guess that for each $p$-block $B$ of maximal defect, the number of height
zero characters in $\irr B$ is equal to the number of height zero characters
in $\irr b$, where $b$ is the Brauer correspondent of $B$.

In fact, it may be unnecessary to limit ourselves to blocks of maximal
defect. Perhaps it is true for {\it every} $p$-block $B$ of $G$ that the
numbers of height zero characters in $\irr B$ and $\irr b$ are equal, where
$b$ is the Brauer correspondent of $B$ with respect to some defect group.
This is precisely the Alperin-McKay conjecture, which appears as
Conjecture~3 of Alperin's paper \ref\alperin. This conjecture has been
shown to be valid for many families of groups. 

Our Conjecture~B strengthens the Alperin-McKay conjecture in the same
way that Conjecture~A strengthens the McKay conjecture. To state it, we
need to define the integer $M_k(B)$, where $B$ is a $p$-block of $G$ and
$k$ is an integer not divisible by $p$. We write $M_k(B)$ to denote the
number height zero characters in $\irr B$ for which the $p'$-part
of the degree is congruent modulo $p$ to $\pm k$. Note that if we hold $k$
fixed and sum $M_k(B)$ over all $p$-blocks $B$ of $G$ of maximal defect,
we obtain the number $M_k(G)$. Also, it is clear that if $p > 2$ and we sum
$M_k(B)$ for $1 \le k \le (p-1)/2$, we get the total number of height zero
characters in $\irr B$.

At this point, one might guess that it is always true that $M_k(B) = M_k(b)$,
where $k$ is any integer not divisible by $p$ and $b$ is the Brauer
correspondent of $B$ with respect to some defect group $D$. To see that
this cannot be correct, however, consider the situation where $G$ is
$p$-solvable. In this case, there is a subgroup $G_0$ of $G$ and a block
$B_0$ of $G_0$ having defect group $D$ and such that the members of
$\irr B$ are exactly the induced characters $\psi^G$, where $\psi$ runs over
$\irr{B_0}$. Now let $N = \norm GD$, write $N_0 = N \cap G_0$ and let $b_0$
be the Brauer correspondent of $B_0$. All of this can be done so that
$D \in \syl p{G_0}$ and the members of $\irr b$ are exactly the induced
characters $\xi^N$, where $\xi$ runs over $\irr{b_0}$. Note that in this
situation, the height zero characters in $\irr B$ are induced from the
height zero characters in $\irr{B_0}$ and similarly, the height zero
characters in $\irr b$ are induced from the height zero characters in
$\irr{b_0}$. 

Suppose now that we knew that $M_k(B_0) = M_k(b_0)$, for all integers $k$
not divisible by $p$. (And in fact, this can be proved using deep facts from
the representation theory of $p$-solvable groups.) Since the height zero
characters of $B$ and $b$ are obtained by induction from $G_0$ and $N_0$,
respectively, it follows that $M_{rk}(B) = M_{sk}(b)$, where $r$ and $s$ are
respectively the $p'$-parts of the indices $|G:G_0|$ and $|N:N_0|$. But in
general, $r$ and $s$ are not congruent modulo $p$, and so we cannot
expect that $M_k(B) = M_k(b)$ for all choices of $k$.

In this situation, let $c$ be the $p'$-part of $|G:N|$. Since
$D \in \syl p{G_0}$, we know by Sylow's theorem that
$|G_0:N_0| \equiv 1$ mod $p$, and thus
$$
r \equiv |G:G_0|_{p'} \equiv
|G:N_0|_{p'} = |G:N|_{p'}|N:N_0|_{p'} \equiv cs {\rm ~~mod~p} \,.
$$
We now have $M_{sk}(b) = M_{rk}(B) = M_{csk}(B)$ for all integers $k$ that
are not divisible by $p$, or equivalently, $M_{ck}(B) = M_k(b)$ for all such
integers $k$. We conjecture that this fact about $p$-solvable groups holds
in general.

\nonumproclaim{Conjecture B} Let $B$ be a $p$\/{\rm -}\/block of an arbitrary finite
group $G$ and suppose that $b$ is the Brauer correspondent of $B$ with
respect to some defect group $D$. Then for each integer $k$ not divisible
by $p${\rm ,} we have $M_{ck}(B) = M_k(b)${\rm ,} where $c = |G:\norm GD|_{p'}$.
\endproclaim

Observe that Conjecture~B implies Conjecture~A. To see why this is so,
recall that $M_k(G)$ is the sum of the quantities $M_k(B)$ as $B$ runs over
all $p$-blocks of $G$ of maximal defect. For each such block, however,
Conjecture~B asserts that $M_k(B) = M_k(b)$, where $b$ is the Brauer
correspondent of $B$ with respect to some fixed Sylow $p$-subgroup $P$ 
of $G$. (This is because by Sylow's theorem, the constant $c$ that appears
in Conjecture~B is congruent to $1$ modulo $p$ in this case.) Conjecture~A
then follows since the sum of the quantities $M_k(b)$ is exactly $M_k(N)$,
where $N = \norm GP$.

As we have already indicated, Conjecture~B is true for $p$-solvable
groups. Also, Fong's paper \ref\fong\ establishes this conjecture for all
symmetric groups (for all primes). There is one other situation where we
know that our Conjecture~B is valid.

\nonumproclaim{(2.1) Theorem} Conjecture~{\rm B} holds for all $p$\/{\rm -}\/blocks that
have cyclic defect groups.
\endproclaim

Since Conjecture~B implies Conjecture~A, it follows that Conjecture~A is
valid for all groups $G$ having a cyclic Sylow $p$-subgroup. (This was
mentioned in Section~1.)

To establish Theorem~2.1, we need to appeal to the deep theory of blocks
with cyclic defect groups that was developed by Dade. A consequence of
this theory in \ref\dade\ is the following useful fact.

\nonumproclaim{(2.2) Theorem} Suppose that $B$ is a $p$\/{\rm -}\/block of $G$ with
cyclic defect group $D = \gen x$ and let $b$ be the Brauer correspondent
of $B$ with respect to $D$. Then for each character $\chi \in \irr B$ there is
a sign $\epsilon_\chi = \pm1$ and a character $\tilde\chi \in \irr b$ such that
$$
\chi(xy) = \epsilon_\chi \tilde\chi(xy)
$$
for all $p$\/{\rm -}\/regular elements $y \in \cent Gx$. Furthermore{\rm ,} the map
$\chi \mapsto \tilde\chi$ defines a bijection from $\irr B$ onto $\irr b$.
\endproclaim

\demo{Sketch of proof} Write $C = \cent GD$ and $N = \norm GD$.
According to Dade's paper \ref\dade, there is a certain uniquely determined
subgroup $E$ with $C \sbs E \sbs N$, where $|E:C| = e$ is a divisor of $p-1$.
(The uniqueness of $E$ depends on the fact that $N/C$ is abelian.) Dade
shows in Theorem~1, Part~1, that the members of $\irr B$ are of two types.
There are $e$ characters of the form $\chi_j$, where $j$ is an integer
$1 \le j \le e$ and there are $(|D|-1)/e$ characters $\chi_\lambda$, where
$\lambda$ is a nonprincipal linear character of $D$. (Our notation, which
differs slightly from Dade's,  is set up so that $\chi_\lambda = \chi_\mu$ if
and only if $\lambda$ and $\mu$ are conjugate in $E$.) Similarly, if we
apply this reasoning to $N$ in place of $G$, with $b$ in place of $B$, we
get the same subgroup $E$, and thus the members of $\irr b$ can be
parametrized as the characters $\psi_j$ with $1 \le j \le e$ and
$\psi_\lambda$, where $\lambda$ is a nonprincipal linear character of $D$
and $\psi_\lambda = \psi_\mu$ if and only if $\lambda$ and $\mu$ are
conjugate in $E$. We can now define the bijection $\chi \mapsto \tilde\chi$
by $\chi_j \mapsto \psi_j$ for $1 \le j \le e$ and $\chi_\lambda \mapsto
\psi_\lambda$ if $\lambda$ is a linear character of $D$.

Dade's Corollary~1.9 gives formulas for the evaluation of $\chi_j(xy)$ and
$\chi_\lambda(xy)$, and of course similar formulas can be used to compute
$\psi_j(xy)$ and $\psi_\lambda(xy)$ by working in $N$ instead of $G$. (Note
that since $D = \gen x$, we must set $i = 0$ in Corollary~1.9.) Examination
of the right sides of the formulas in Corollary~1.9 shows that everything is
determined inside the group $N$ except for a constant $\gamma_0$ (which
is always equal to $1$ by Dade's Equation~1.10) and certain signs
$\epsilon_j$, which depend on the character. (All of the characters
$\chi_\lambda$ use the same sign, $\epsilon_0$.) It follows that $\chi(xy)$
and $\tilde\chi(xy)$ agree, except possibly for a sign depending on $\chi$.
This completes the proof.\hfill\qed
\enddemo

{\it Proof of Theorem} 2.1. Let $B$ be a $p$-block of $G$ with cyclic
defect group $D = \gen x$ and let $b$ be the Brauer correspondent of $B$
with respect to $D$, so that $b$ is a $p$-block of $N = \norm GD$. We will
show that the bijection $\chi \mapsto \tilde\chi$ of Theorem~2.2 maps
the height zero characters in $\irr B$ onto the height zero characters in
$\irr b$. (Actually, it is true that in this case all of the members of $\irr B$
and $\irr b$ have height zero, but we will not need that fact.) Also, we will
show that if $\chi$ and $\tilde\chi$ are height zero characters, then
$\chi(1)_{p'} \equiv \pm c\tilde\chi(1)_{p'}$ mod $p$, where
$c = |G:N|_{p'}$. The result will then follow. 

Let $K$ be a defect class for $B$. In particular, $D$ is a defect group for
$K$, which means that there exists $y \in K$ such that
$D \in \syl p{\cent Gy}$, and thus $y \in C$, where
$C = \cent GD = \cent Gx$. Also, because $K$ is a defect class for $B$, we
know that $y$ is $p$-regular and that $\lambda_B(\hat K) \ne 0$, where
$\lambda_B$ is the ``central homomorphism" corresponding to $B$ and
$\hat K$ is the sum of the elements of $K$ in the appropriate group ring.
(Recall that for every class $L$ of $G$, we have
$\lambda_B(\hat L) =  \omega_\chi(\hat L)^*$, where $\chi$ is any member
of $\irr B$ and $(\ )^*$ is the canonical homomorphism from the ring of
$p$-local integers to its residue class field modulo some fixed maximal ideal
$M$ containing $p$.)

Now let $L = \cl G{xy}$ and note that $D \in \syl p{\cent G{xy}}$, so that
$|L|_p = |G|_p/|D| = |K|_p$. We claim now that
$\lambda_B(\hat L) \ne 0$. To see why this is so, let $\chi \in \irr B$ have
height zero. Then $|K|_p = \chi(1)_p = |L|_p$, and thus $|K|/\chi(1)$ and
$|K|/|L|$ are $p$-local integers. Also, since $\chi(y) \equiv \chi(xy)$,
where we are working modulo the maximal ideal $M$, it follows that 
$$
\omega_\chi(\hat K) =
{\chi(y)|K| \over \chi(1)} \equiv
{\chi(xy)|K| \over \chi(1)} =
{\chi(xy)|L| \over \chi(1)}\,{|K| \over |L|} =
\omega_\chi(\hat L){|K| \over |L|} \,.
$$
We now have
$$
0 \ne \lambda_B(\hat K) = \omega_\chi(\hat K)^* =
\omega_\chi(\hat L)^* \left({|K| \over |L|}\right)^ {\kern -.3em *} =
\lambda_B(\hat L) \left({|K| \over |L|}\right)^ {\kern -.3em *}  \,,
$$
and it follows that $\lambda_B(\hat L) \ne 0$, as claimed.

By Lemma~15.46 of \ref\book, we know that $L \cap C$ is a class of $N$,
and thus $L \cap C = \cl N{xy}$. Since $\cent G{xy} \sbs C \sbs N$, we see
that $\cent G{xy} = \cent N{xy}$, and from this, we compute that
$|L| = |G:\cent G{xy}| = |G:N||N:\cent N{xy}| = |G:N||L \cap C|$.
Observe that since $b^G = B$, we have
$\lambda_B(\hat L) = \lambda_b(\widehat{L \cap C})$, and we write
$\alpha$ to denote this nonzero element of the residue class
field of the $p$-local integers.

Now let $\chi \in \irr B$ be arbitrary and write $\psi = \tilde\chi$, in the
notation of Theorem~2.2. We thus have
$$
\omega_\chi(\hat L)^* = \alpha = \omega_\psi(\widehat{L \cap C})^* \,.
$$
Since $L \cap C$ is a class of $N$ with defect group $D$, we see that
$\chi(1)/|L|$ and $\psi(1)/|L \cap C|$ are $p$-local integers. Also, we
observe that $\chi$ has height zero in $B$ precisely when
$(\chi(1)/|L|)^* \ne 0$, and similarly, $\psi$ has height zero in $b$ if and
only if $(\psi(1)/|L \cap C|)^* \ne 0$.

By Theorem~2.2, we have
$$
{\chi(1) \over |L|} \omega_\chi(\hat L) = \chi(xy) = \pm \psi(xy) =
\pm{\psi(1) \over |L \cap C|} \omega_\psi(\widehat{L \cap C}) \,.
$$
Since $\chi(1)/|L|$ and $\psi(1)/|L \cap C|$ are $p$-local integers, we
deduce that 
$$
\left({\chi(1) \over |L| }\right)^{\kern -.3em *}\!\alpha =
\pm \left({\psi(1) \over |L \cap C|}\right)^ {\kern -.3em *}\!\alpha \,,
$$
and thus since $\alpha \ne 0$, we have
$\chi(1)/|L| \equiv \pm \psi(1)/|L \cap C|$. In particular, $\chi$ has height
zero if and only if $\psi$ has height zero. If we multiply both sides by
$|L|_{p'} = |G:N|_{p'}|L \cap C|_{p'} = c|L \cap C|_{p'}$, we obtain
$$
{\chi(1) \over |L|_p} \equiv \pm c {\psi(1) \over |L \cap C|_p}
$$
and since both sides of this congruence are rational integers, these
numbers are actually congruent modulo $p$. In the case where $\chi$ and
$\psi$ have height zero, the integers $\chi(1)/|L|_p$ and
$\psi(1)|L \cap C|_p$ are exactly the $p'$-parts of the degrees of $\chi$
and $\psi$, and so the result follows.\hfill\qed
\vglue12pt

Our Conjecture~B is related to some conjectures and results of M.~Brou\'e
in \ref\broue. As usual, suppose that $B$ is a $p$-block of $G$ and that $b$
is its Brauer correspondent with respect to the defect group $D$. Brou\'e
conjectures that in the case where $D$ is abelian, there exists a ``perfect
isometry" between $B$ and $b$, and he proves this in the case where $D$
is cyclic. (It is known that a perfect isometry need not exist in the case
where $D$ is nonabelian.) A perfect isometry, if it exists, would imply the
existence of a certain bijection $\chi \mapsto \tilde\chi$ from $\irr B$ onto
$\irr b$. In addition, Brou\'e shows that if a perfect isometry exists, there
would be some constant $c$, depending on the block $B$, such that
$\chi(1) \equiv \pm c \tilde\chi(1)$ mod $p$ for all $\chi \in \irr B$. Of
course, it follows in this case that for each integer $k$ not divisible by $p$,
we would have (in our notation) $M_{ck}(B) = M_k(b)$. Furthermore,
Brou\'e shows that his constant $c$ is equal to $1$ if $B$ is the principal
block of $G$, but he does not evaluate the constant in other cases. (Note
that according to Conjecture~B, this constant should be $1$ for every
block of maximal defect, and not just for the principal block.) We mention
that if $G$ has an abelian self-centralizing Sylow $p$-subgroup, then the
principal block is the only $p$-block of maximal defect, and in that case,
Brou\'e's perfect isometry conjecture would imply our Conjecture~A.

\section{Field automorphisms}

There are other directions in which the McKay conjecture might be
extended. Suppose, for example, that $G$ is a group for which the McKay
conjecture holds in the strong sense that there is a {\it canonical} bijection
$\chi \mapsto \tilde\chi$ from $\irrpp G$ onto $\irrpp N$. (Here
$N = \norm GP$, where $P \in \syl pG$, and we are using the notation
$\irrpp X$ to denote the subset of $\irr X$ consisting of characters of
$p'$-degree.) In this case, we see that if $\sigma$ is any automorphism of
the cyclotomic field $\Q_{|G|}$, then $(\tilde\chi)^\sigma =
\widetilde{\chi^\sigma}$, and thus in particular, the sets $\irrpp G$ and
$\irrpp N$ would have equal numbers of $\sigma$-fixed members. 

If $G$ is (solvable) of odd order, then there is such a canonical bijection,
and the numbers of $\sigma$-fixed characters in $\irrpp G$ and $\irrpp N$
are equal for all choices of the field automorphism $\sigma$. But this fails
for solvable groups in general. (For example, if $G = {\rm GL}(2,3)$ and $p = 3$,
then all members of $\irrpp N$ are rational valued, but the same is not true
for $\irrpp G$.) If we impose some conditions on the field automorphism
$\sigma$, however, then the equality of the numbers of $\sigma$-fixed
characters is known to hold for all $p$-solvable groups. (See Corollary~C of
\ref{\isnav}.) We conjecture that under these conditions on $\sigma$,
equality holds for all groups.

\nonumproclaim{Conjecture C}  Let $G$ be an arbitrary finite group and fix a
prime~$p$. Let $\sigma$ be an automorphism of the cyclotomic field
$\Q_{|G|}$ and assume that $\sigma$ has $p$\/{\rm -}\/power order and that it fixes
all $p'$\/{\rm -}\/roots of unity in $\Q_{|G|}$. Then $\sigma$ fixes equal numbers of
characters in $\irrpp G$ and $\irrpp N$, where $N = \norm GP$ and $P \in
\syl pG$.
\endproclaim

Of course, if we take $\sigma$ to be the identity automorphism, we
recover the McKay conjecture from Conjecture~C. Another consequence
of the conjecture is that the character table of a group $G$ determines the
exponent of the abelian group $P/P'$, where, $P \in \syl pG$. To see why
this is true, let $N = \norm GP$ and fix a positive integer $n$. Let
$\sigma_n$ be the unique automorphism of $\Q_{|G|}$ that fixes all
$p'$-roots of unity and maps every $p$-power root of unity $\epsilon$ to
$\epsilon^{p^n+1}$. Then $\sigma_n$ has $p$-power order and it fixes roots
of unity of order $p^n$ but not those of order $p^{n+1}$ or higher. It is not
hard to see from this that a necessary and sufficient condition for
$\sigma_n$ to fix every member of $\irrpp N$ is that $P/P'$ has exponent at
most $n$. If Conjecture~C is true, therefore, it follows that $P/P'$ has
exponent at most $n$ if and only if $\sigma_n$ fixes every member of
$\irrpp G$, and thus the exponent of $P/P'$ is determined from the
character table of $G$, as claimed.

We also propose a block version of Conjecture~C that generalizes the
Alperin-McKay conjecture. To state it, we observe that if $\sigma$ is a field
automorphism fixing $p'$-roots of unity and $\chi \in \irr G$, then $\chi$
and $\chi^\sigma$ necessarily belong to the same $p$-block of $G$. (This is
because these characters agree on all $p'$-elements of $G$.) Such a field
automorphism, therefore, permutes the set of height zero characters in
$\irr B$ for each $p$-block $B$ of $G$.

\nonumproclaim{Conjecture D}  Let $B$ be a $p$-block for an arbitrary finite
group $G$ and suppose that $b$ is the Brauer correspondent of $B$ with
respect to some defect group. Let $\sigma$ be an automorphism of the
cyclotomic field $\Q_{|G|}$ and assume that $\sigma$ has $p$\/{\rm -}\/power order
and that it fixes all $p'$\/{\rm -}\/roots of unity in $\Q_{|G|}$. Then $\sigma$ fixes
equal numbers of height zero characters in $\irr B$ and $\irr b$.
\endproclaim

By Theorem~G of \ref\isnav, Conjecture~D is known to hold for
$p$-solvable groups. We present a proof of the conjecture in the case
where the defect group of $B$ is cyclic.

\nonumproclaim{(3.1) Theorem} Conjecture~$D$ is valid for blocks with cyclic
defect group.
\endproclaim

{\it Sketch of proof}. We suppose that the defect group $D$ of $B$ is
cyclic and that $b$ is the Brauer correspondent of $B$ with respect to $D$.
As we observed in the proof of the Theorem~2.2, the members of $\irr B$
are of two types. There are $e$ characters $\chi_j$, where $1 \le j \le e$
and there are $(|D| - 1)/e$ different characters of the form
$\chi_\lambda$, where $\lambda$ is a nonprincipal linear character of $D$
and $\chi_\lambda = \chi_\mu$ if and only if $\lambda$ and $\mu$ lie in the
same $E$-orbit. Also, there is a similar parametrization of the members of
$\irr b$.

We will see that all of the characters $\chi_j$ and $\psi_j$ are fixed by
$\sigma$ and that $\chi_\lambda$ and $\psi_\lambda$ are fixed by $\sigma$
if and only if the linear character $\lambda$ is $\sigma$-fixed. The result
will then follow.

Suppose that $\chi \in \irr B$ and that $g \in G$ is arbitrary. If $g$
is $p$-regular, we know that $\chi(g)$ is $p$-rational, and thus
$\chi(g) = \chi(g)^\sigma$. Also, if the $p$-part of $g$ is not conjugate to
an element of $D$, then $\chi(g) = 0 = \chi(g)^\sigma$. It follows that to
determine whether or not $\chi$ is $\sigma$-fixed, it suffices to consider
only the values $\chi(xy)$, where $1 \ne x \in D$ and $y$ is a $p$-regular
element in $\cent Gx$. We are thus exactly in the situation where
Corollary~1.9 of \ref\dade\ applies.

It is immediate from Dade's Corollary~1.9 that if $\chi = \chi_j$ with
$1 \le j \le e$, then $\chi(xy)$ is $p$-rational, and it follows that all of the
characters $\chi_j$ are $\sigma$-fixed, as claimed. Also from Corollary~1.9,
we see that $\chi_\lambda(xy)^\sigma = \chi_{\lambda^\sigma}(xy)$. 

If $\lambda$ is $\sigma$-fixed, it is now immediate that $\chi_\lambda$ is
$\sigma$-fixed. Conversely, if $\chi_\lambda$ is $\sigma$-fixed and we
write $\mu = \lambda^\sigma$, then we have $\chi_\lambda = \chi_\mu$,
and thus the $E$-orbit containing $\lambda$ is invariant under $\sigma$.
Since $\sigma$ has $p$-power order, however, and the $E$-orbit
containing $\lambda$ has size $e$, which is not divisible by $p$, it follows
that every member of this $E$-orbit is $\sigma$-fixed, and in particular
$\lambda$ is $\sigma$-fixed, as desired.

We have now shown that as claimed, the $\sigma$-fixed members of $\irr B$
are exactly the characters $\chi_j$ with $1 \le j \le e$ and the characters
$\chi_\lambda$, where $\lambda$ is $\sigma$-fixed. Exactly the same
reasoning applies to the block $b$, and the proof is complete.\hfill\qed
 \vglue9pt

Of course, one could combine our Conjectures~A and~C. Perhaps, for
example, it is true that for each integer $k$ not divisible by $p$ and each
appropriate field automorphism $\sigma$, the groups $G$ and
$N = \norm GP$ always have equal numbers of $\sigma$-fixed characters
with degrees congruent modulo $p$ to $\pm k$. Similarly, one could
combine our Conjectures~B and~D, but we will refrain from stating these
composite conjectures formally.
\bigskip\bigskip

\references

alperin
 \name{J.\ L.\ Alperin}, The main problem of block theory,  {\it Proc.\ 
of the Conference on Finite Groups\/} (Univ. of Utah, Park City, Utah,
1975),
 341--356, Academic Press, New York, 1976.

2
 \name{M.\ Brou\'e},  Isom\'etries parfaites, types de blocs, cat\'egories
d\'eriv\'ees, {\it Ast{\rm \'{\it e}}risque\/} {\bf 181--182} (1990), 61--92.

3
 \name{E.\ C.\ Dade},  Blocks with cyclic defect groups, {\it Ann.\ of
Math\/}.\  
{\bf 84 } (1966), 20--48.

4
\name{P.\ Fong},  The Isaacs-Navarro conjecture for symmetric groups, {\it J. Algebra}, to appear.

5
 \name{P.\ Fong} and \name{B.\ Srinivasan}, The blocks of finite general linear and
unitary groups, {\it Invent. Math}.\/   {\bf 69} (1982), 109--153.

6
  \name{I.\ M.\ Isaacs},  Characters of solvable and symplectic groups,
{\it Amer.\ J.\  Math\/}.\  {\bf 95} (1973), 594--635.

7
  \bibline, {\it Character Theory of Finite Groups\/}, Dover Publ.\ Inc.,
New York, 1994.

8
  \name{I.\ M.\ Isaacs} and \name{G.\ Navarro}, Characters of $p'$-degree of
$p$-solvable groups, {\it J.\ Algebra\/}, {\bf 246} (2001), 394--413.

9
   \name{J.\ McKay}, A new invariant for simple groups, {\it Notices
Amer.\
Math.\ Soc\/}.\  {\bf 18} (1971), 397.

10   \bibline, Irreducible representations of odd degree, {\it J.\
Algebra\/}
 {\bf 20} (1972), 416--418.
 
11   \name{R.\ A.\ Wilson},  The McKay conjecture is true for the sporadic
simple groups, {\it J.\  Algebra\/} {\bf 207} (1998), 294--305.

\endreferences
\bye